\documentclass[11pt]{amsart2000}
\usepackage{amssymb}
\usepackage{enumerate}
\newtheorem{theorem}{Theorem}[section]
\newtheorem{lemma}[theorem]{Lemma}
\newtheorem{corollary}[theorem]{Corollary}
\newtheorem{problem}[theorem]{Problem}
\newcommand{\easytosee}{$\square$}
\theoremstyle{definition}
\newtheorem{definition}[theorem]{Definition}
\newtheorem{example}[theorem]{Example}
\theoremstyle{remark}
\newtheorem{remark}[theorem]{Remark}

\newtheorem{qlemma}{Lemma}

\begin{document}
\setlength{\unitlength}{0.01in}
\linethickness{0.01in}
\begin{center}
\begin{picture}(474,66)(0,0)
\multiput(0,66)(1,0){40}{\line(0,-1){24}}
\multiput(43,65)(1,-1){24}{\line(0,-1){40}}
\multiput(1,39)(1,-1){40}{\line(1,0){24}}
\multiput(70,2)(1,1){24}{\line(0,1){40}}
\multiput(72,0)(1,1){24}{\line(1,0){40}}
\multiput(97,66)(1,0){40}{\line(0,-1){40}}
\put(143,66){\makebox(0,0)[tl]{\footnotesize Proceedings of the Ninth Prague Topological Symposium}}
\put(143,50){\makebox(0,0)[tl]{\footnotesize Contributed papers from the symposium held in}}
\put(143,34){\makebox(0,0)[tl]{\footnotesize Prague, Czech Republic, August 19--25, 2001}}
\end{picture}
\end{center}
\vspace{0.25in}
\setcounter{page}{93}
\title[Characterizing continuity]{Characterizing continuity by preserving
compactness and connectedness}
\author[J. Gerlits]{J\'anos Gerlits}
\address{Alfr\'ed R\'enyi Institute of Mathematics\\
P.O.Box 127, 1364 Budapest, Hungary}
\email{gerlits@renyi.hu}
\author[I. Juh\'asz]{Istv\'an Juh\'asz}
\address{Alfr\'ed R\'enyi Institute of Mathematics\\
P.O.Box 127, 1364 Budapest, Hungary}
\email{juhasz@renyi.hu}
\author[L. Soukup]{Lajos Soukup}
\address{Alfr\'ed R\'enyi Institute of Mathematics\\
P.O.Box 127, 1364 Budapest, Hungary}
\email{soukup@renyi.hu}
\author[Z. Szentmikl\'ossy]{Zolt\'an Szentmikl\'ossy}
\address{E\"otv\"os Lor\'ant University, Department of Analysis,
1117 Budapest, P\'azm\'any P\'eter s\'et\'any 1/A, Hungary}
\email{zoli@renyi.hu}
\thanks{The second author was an invited speaker at the Ninth Prague Topological Symposium.}
\thanks{Research supported by Hungarian Foundation for Scientific
Research, grant No.25745.}
\thanks{This article has been submitted for publication to
\textit{Fundamenta Mathematicae}.}
\thanks{J\'anos Gerlits, Istv\'an Juh\'asz, Lajos Soukup and Zolt\'an 
Szentmikl\'ossy,
{\em Characterizing continuity by preserving compactness and 
connectedness},
Proceedings of the Ninth Prague Topological Symposium, (Prague, 2001),
pp.~93--118, Topology Atlas, Toronto, 2002}
\begin{abstract}
Let us call a function $f$ from a space $X$ into a space
$Y$ {\it preserving} if the image of every compact subspace of $X$
is compact in $Y$ and the image of every connected subspace of $X$
is connected in $Y$. By elementary theorems a continuous function
is always preserving. Evelyn R. McMillan \cite{McM} proved in 1970
that if $X$ is Hausdorff, locally connected and Fr\`echet, $Y$ is
Hausdorff, then the converse is also true: any preserving function
$f:X\to Y$ is continuous.
The main result of this paper is that if $X$ is any product of
connected linearly ordered spaces (e.g.\ if $X={\mathbb R}^\kappa$) and
$f:X \to Y$ is a preserving function into a
regular space $Y$, then $f$ is continuous.
\end{abstract}
\maketitle

Let us call a function $f$ from a space $X$ into a space $Y$ {\it
preserving} if the image of every compact subspace of $X$ is
compact in $Y$ and the image of every connected subspace of $X$ is
connected in $Y$. By elementary theorems a continuous function is
always preserving. Quite a few authors noticed---mostly
independently from each other---that the converse is also true for
real functions: a preserving function $f:{\mathbb R} \to
{\mathbb R}$ is continuous. (The first paper we know of is
\cite {Rowe}
from 1926!) Whyburn proved \cite{Wh} that a preserving function
from a space $X$
into a Hausdorff space is always continuous at a first countability and
local connectivity point of $X$.

Evelyn R. McMillan \cite{McM} proved in
1970 that if $X$ is Hausdorff, locally connected and Fr\`echet,
moreover $Y$ is Hausdorff, then any preserving function $f:X\to Y$ is
continuous. This is quite a significant and deep result that is
surprisingly little known.

We shall use the notation $Pr(X,T_i)$
($i=1,2,3$ or $3{\frac{1}{2}}$) to denote the following statement:
Every preserving function from the topological space $X$ into
any $T_i$ space is continuous.

The organization of the paper is as follows: In \S $1$ we give some basic
definitions and then treat our results that are closely connected to
McMillan's theorem. 
\S $2$ treats several important technical theorems that enable us to
conclude that certain preserving functions are continuous.
In \S $3$ we prove that certain product spaces $X$ satisfy $Pr(X,T_3)$; in
particular, any preserving function from a product of connected linearly
ordered spaces into a regular space is continuous. 
In \S $4$ we discuss some results concerning the continuity of preserving
functions defined on compact or sequential spaces. 
Finally, \S $5$ treats the relation $Pr(X,T_1)$.

Our terminology is standard.
Undefined terms can be found in \cite{Eng} or in \cite{Juh}.

\section{Around McMillan's theorem}

Our first theorem implies that (at least among Tychonoff spaces) local
connectivity of $X$ is a necessary condition for 
$Pr(X,T_{3{\frac{1}{2}}})$. 
For metric spaces this result was proved by Klee and Utz\cite{Kl}.

\begin{theorem}
\label{1.1}
If the Tychonoff space $X$ is not locally connected at a point $p \in X$,
then there exists a preserving function $f$ from $X$ into the interval
$[0,1]$ which is not continuous at $p$.
\end{theorem}

\begin{proof} 
Suppose $X$ is not locally connected at the point $p$, then there is an
open neighbourhood $U$ of $p$ such that if $K$ denotes the component of
$p$ in $U$, then $K$ is not a neighbourhood of $p$. As $X$ is Tychonoff,
there exists an open neighbourhood $V\subset U$ of $p$ and a continuous
function $\overline f:X \to [0,1]$ such that $\overline f$ is identically
$0$ on $V$ and identically $1$ outside of $U$. Select another continuous
function $g:X \to [0,1]$ such that $g(p)=1$ and $g$ is identically $0$
outside $V$. Put now
\begin{displaymath}
f(x)=
\left\{
\begin{array}{ll}
\overline{f}(x) + g(x)&\text{ if $x\in K$;}\\
\overline {f}(x) & \text{otherwise.}
\end{array}
\right.
\end{displaymath}

($f:X \to [0,1]$ because $\overline f[V]=\{0\}$ and $g[X-V]=\{0\}$.) 
The function $f$ is not continuous at $p$ because $f(p)=1$ and in every
neighbourhood of $p$ there is a point from $V - K$ mapped to $0$.

On the other hand, we claim that $f$ is preserving. Indeed, let $C \subset
X$ be compact. The restriction of $f$ to the closed set $F=(X-V)\cup K$ is
evidently continuous (here $f=\overline f + g$), hence $C'=f(C\cap F)$ is
compact. But $f$ is $0$ on $V-K$ and so $f(C)$ is either $C'$ or $C' \cup
\{0\}$, thus $f(C)$ is clearly compact.

Let now $C$ be any connected subset of $X$. The function $f$ is continuous
on $X-K$ (here $f=\overline f$), hence $C\cap K \not = \emptyset$ can be
assumed. Similarly, $f$ is continuous on $X-V$ ($f(x)=\overline f(x)$ for
$x\in X-V$), hence we can suppose that $C$ also meets $V$. If $C \subset
U$, then necessarily $C \subset K$ because $C$ is connected and meets the
component $K$of $U$. As $f$ is continuous on $K$, the image $f(C)$ is then
connected.

Thus it remains to check only the case in which $C$ meets both $V$ and
$X-U$. In this case $\overline f(C)$ is the whole interval $[0,1]$. But
$f$ and $\overline f$ are equal at each point where $\overline f$ does not
vanish, thus $f(C)$ contains the interval $(0,1]$ and so it is connected.
\end{proof}

It is not a coincidence that the target space in Theorem \ref{1.1} is
the interval $[0,1]$, because of the following result:

\begin{lemma}
\label{1.2} 
Suppose $f: X \to Y$ is a preserving function into a Tychonoff space $Y$
and $f$ is not continuous at the point $p \in X$. Then there exists a
preserving function $h:X \to [0,1]$ which is also not continuous at $p$.
\end{lemma}

\begin{proof} 
Since $f$ is not continuous at $p$, there exists a closed set $F\subset Y$
such that $f(p)\not \in F$ but $p$ is an accumulation point of
$f^{-1}(F)$. 
Choose a continuous function $g:Y\to [0,1]$ such that $g(f(p))=0$ and $g$ 
is identically $1$ on $F$. Then the composite function $h(x)=f(g(x))$ has
the stated properties.
\end{proof}

The following Lemmas will be often used in the sequel. 
They all state simple properties of preserving functions.

\begin{lemma}
\label{1.3}If $f:X \to Y$ is a
compactness preserving function, $Y$ is Hausdorff, $M
\subset X$ with $\overline M$ compact and $f(M)$
infinite then for every accumulation point $y$ of $f(M)$
there is an accumulation point $x$ of $M$ such that
$f(x) = y$, i. e. $f(M)' \subset f(M')$ .
\end{lemma}

\begin{proof} 
Let $N = M - f^{-1}(y)$ then $f(N) = f(M) - \{y\}$ and so we have 
$y\in \overline {f(N)}- f(N)$. But $f(\overline {N})$ is also compact,
hence closed in $Y$ and so $y\in f(\overline {N}) - f(N)$ as well.
Thus there is an $x\in \overline N - N$ such that $f(x)=y$ and then $x$ is
as required.
\end{proof}

We shall often use the following immediate consequence of this lemma:

\begin{qlemma}[E.~R.~McMillan \cite{McM}]
\label{1.3'}
If $f:X \to Y$ is a compactness preserving function, $Y$ is Hausdorff,
$\{x_n : n<\omega \} \subset X$ converges to $x\in X$ then either 
${\{ f(x_n): n<\omega \}}$, converges to $f(x)$ or the set
$\{f(x_n): n<\omega \}$ is finite
\easytosee
\end{qlemma}

Actually, to prove Lemma \ref{1.3'} we do not need the full force of the
assumption that $f$ is compactness preserving. It suffices to assume that
the image of a convergent sequence together with its limit is compact, in
other words: the image of a topological copy of $\omega +1$ is compact.
For almost all of our results given below only this very restricted
special case of compactness preservation is needed.

\begin{lemma}[\cite{Per}]
\label{1.4} 
If $f:X \to Y$ preserves connectedness, $Y$ is a $T_1$-space and $C\subset
X$ is a connected set, then $f(\overline C) \subset \overline {f(C)}$.
\end{lemma}

\begin{proof}
If $x\in \overline C$ then $C\cup \{x\}$ is connected. 
Thus $f(C\cup \{x\})= f(C)\cup \{f(x)\}$ is also connected and hence
$f(x) \in \overline {f(C)}$.
\end{proof}

The next lemma will also play a crucial role in some theorems of the
paper. A weaker form of it appears in \cite{McM}.

\begin{definition}
\label{1.5} 
We shall say that $f:X \to Y$ is {\it locally constant} at the point 
$x \in X$ if there is a neighbourhood $U$ of $x$ such that $f$ is constant
on $U$.
\end{definition}

\begin{lemma}
\label{1.6}
Let $f$ be a connectivity preserving function from a locally connected
space $X$ into a $T_1$-space $Y$. If $F\subset Y$ is closed and $p \in
\overline {f^{-1}(F)} - f^{-1}(F)$ then $p$ is also in the closure of
the set
\begin{displaymath}
\{ x \in f^{-1}(F): f \hbox { is not locally constant at
} x \}.
\end{displaymath}
\end{lemma}

\begin{proof} 
Let $G$ be a connected open neighbourhood of $p$ and $C$ be
a component of the non-empty subspace $G \cap f^{-1}(F)$. Then $C$ has a
boundary point $x$ in the connected subspace $G$ because $\emptyset \not =
C \not = G$. Clearly, $f(x) \in F$ by Lemma \ref{1.4}. If $V\subset G$ is
any connected neighbourhood of $x$ then $V \cup C$ is connected and $V-C
\not = \emptyset$ because $x$ is a boundary point of $C$ hence $V$ is not
contained in $f^{-1}(F)$, so $f$ is not locally constant at $x$.
\end{proof}

\begin{lemma}
\label{1.7} 
Let $f:X \to Y$ be a connectivity preserving function into the $T_1$-space
$Y$.
Suppose that $X$ is locally connected at the point $p \in X$ and $f$ is
not locally constant at $p$. 
Then $f(U)\cap V$ is infinite for every neighbourhood $U$ of $p$ and for
every neighbourhood $V$ of $f(p)$.
\end{lemma}

\begin{proof}
Choose any connected neighbourhood $U$ of $x$; then $f(U)$ is connected
and has at least two points. Thus if $V$ is any open subset of $Y$ 
containing $f(p)$ then $f(U) \cap V$ can not be finite because otherwise
$f(p)$ would be an isolated point of the non-singleton connected set
$f(U)$.
\end{proof}

The following proof of McMillan's theorem is based upon the same ideas as
her original proof, although, we think, it is much simpler. 
We include it here mainly to make the paper self-contained.

\begin{theorem}[E.~R.~McMillan \cite{McM}]
\label{1.8}
If $X$ is a locally connected and Fr\`echet Hausdorff space,
then $Pr(X,T_2)$ holds.
\end{theorem}

\begin{proof} Assume $Y$ is $T_2$ and $f: X \to Y$ is
preserving but not continuous at the point $p \in X$. Then by
Lemma \ref{1.3'} there is a sequence $x_n \to p$ such that $f(x_n) = y
\not = f(p)$ for all $n<\omega$ . Using Lemma \ref{1.6}
with $F = \{y\}$ we can also
assume that $f$ is not locally constant at the points $x_n$.

As $Y$ is $T_2$, there is an open set $V \subset Y$ such that $y
\in V$ but $f(p) \not \in \overline V$. By Lemma \ref{1.7} the image of
every neighbourhood of each point $x_n$ contains infinitely many
points (different from $y$) from $V$.

Now we select recursively sequences $\{x^n_k : k<\omega \}$
converging to $x_n$ for all $n<\omega$. Suppose $n<\omega$ and
the points $x^m_k$ are already defined for $m<n$ and $k<\omega$ so
that $f(x^m_k)\not = y$. Then $x_n$ is in the closure
of the set
$$f^{-1}(V-(\{f(x^m_k) : m,k<n \}\cup \{y\})),$$
hence, as $X$ is Fr\'echet , the new
sequence $\{x^n_k : k<\omega \}$ converging to $x_n$ can be chosen
from this set.

Since the sequence $\{x^n_k : k<\omega \}$ converges to $x_n$ and $\{x_n
:n<\omega \}$ converges to the (Fr\`echet) point $p$, there is also
a sequence $\{x^{n_l}_{k_l} : l<\omega \}$ converging to $p$. But then
the sequence $\{n_l :
l<\omega \}$ must tend to infinity so, by passing to a subsequence
if necessary, we can assume that $n_{l+1}>\max (n_l,k_l)$ for all
$l<\omega$. However, the sequence $\{f(x^{n_l}_{k_l}) : l<\omega \}$
does not converge to $f(p)$ because
$f(p)\not \in \overline {\{f(x^{n_l}_{k_l}):
l<\omega \}} \subset \overline V $,
while the points $f(x^{n_l}_{k_l})$ are all
distinct, contradicting Lemma \ref{1.3'}. \end{proof}

We could prove the following local version of McMillan's theorem:

\begin{theorem}
\label{1.9} 
If $X$ is a locally connected Hausdorff space, $p$ is a Fr\`echet point of
$X$ and $f$ is a preserving function from $X$ into a Tychonoff space $Y$,
then $f$ is continuous at $p$.
\end{theorem}

\begin{proof}
By Lemma \ref{1.2} it suffices to prove this in the case
when $Y$ is the interval $[0,1]$.
Assume, indirectly, that $f$ is not continuous at $p$
then, since $p$ is a Fr\'echet point and by Lemma \ref{1.6}, we can choose a
sequence $x_n \to p$ and a $y \in [0,1]$ with $y\not= f(p)$
such that $f(x_n) = y$ and $f$ is not locally constant at $x_n$
for all $n < \omega$.

For each $n$ choose a neighbourhood $U_n$ of $x_n$ with $p \not \in
\overline U_n$ and put $A_n = \{x\in U_n : 0<|f(x)-y|<1/n \}$. For any
connected neighbourhood $W$ of $x_n$ its image $f(W)$ is a non-singleton
interval containing $y$, hence the local connectivity of $X$ implies that
$x_n \in \overline A_n$ for all $n<\omega$ and so $p$ belongs to the
closure of $\bigcup \{A_n : n<\omega \}$. As $p$ is a Fr\`echet point,
there is a sequence $z_k \in A_{n_k}$ converging to $p$. But $n_k$
necessarily tends to infinity because $p \not \in \overline A_n \subset
\overline U_n$ for each $n < \omega$, hence $f(z_k) \to y \not = f(p)$,
contradicting Lemma \ref{1.3'}. Indeed, the set $\{f(z_k) : k< \omega \}$
is infinite because $f(z_k)\not = y$. 
\end{proof}

Theorem \ref{1.9} is not a full local version of Theorem \ref{1.8}
because local connectivity is assumed in it globally for X.
This leads to the following natural question:

\begin{problem}
\label{1.10} Let $X$ be a Hausdorff space,
$f$ be a preserving function from $X$ into a Tychonoff space $Y$
and let $X$ be locally connected and Fr\`echet at the point $p\in
X$. Is it true then that $f$ is continuous at $p$?
\end{problem}

We do not know the answer to this problem, however
we can prove some partial affirmative results.

\begin{definition}[\cite{Arh}]
\label{1.11} A point $x$ of
a space $X$ is called an $(\alpha_4)$ point if for any sequence
$\{A_n :n<\omega \}$ of countably infinite sets with $A_n \to x$
for each $n<\omega$ there is a countably infinite set $B \to x$
such that $\{n<\omega: A_n \cap B \not = \emptyset \}$ is
infinite.

An $(\alpha_4)$ and Fr\`echet point will be called an
$(\alpha_4)$-F point in $X$.

\end{definition}

\begin{theorem}
\label{1.12} Let $f$ be a preserving
function from a topological space $X$ into a Hausdorff space $Y$
and let $p$ be a point of local connectivity and an $(\alpha_4)$-F
point in $X$. Then $f$ is continuous at $p$.
\end{theorem}

\begin{proof} Assume not. Then by the Lemma \ref{1.3'}
there is a point $y\in Y$ such that $y\not = f(p)$ but $p$ is in
the closure of $f^{-1}(y)$. Choose an open neighbourhood $V$
of $y$ in $Y$ with $f(p) \not \in \overline V$. By Lemma \ref{1.7}
and Lemma \ref{1.3'} we can recursively choose
pairwise distinct points
$y_n\in V$ such that $p$ is in the closure of $f^{-1}(y_n)$
for all $n\in \omega$. As the point $p$ is an $(\alpha_4)$-F
point in $X$, there is a ``diagonal'' sequence $\{x_m:m\in M \}$
converging to $p$, where $f(x_m)=y_m$ and $M\subset \omega$ is
infinite, contradicting Lemma \ref{1.3'}.\end{proof}

The next result yields a different kind of
partial answer to Problem \ref{1.10}:

\begin{theorem}
\label{1.13} Let $f$ be a preserving
function from a topological space $X$ into a Tychonoff space $Y$
and let $p$ be a Fr\` echet point of local connectivity of $X$
with character $\le 2^{\omega}$. Then $f$ is continuous at $p$.
\end{theorem}

\begin{proof} Assume not, $f$ is discontinuous at the
point $p \in X$. By Lemmas \ref{1.2} and \ref{1.3'}
we can suppose that $Y=[0,1]$,
$f(p)=0$ and every neighbourhood of $p$ is mapped onto the whole
interval $[0,1]$ .
Let $\mathcal U$ be a neighbourhood-base of $p$ of size $\le
2^{\omega}$ and choose for each $U \in \mathcal U$ a point $x_U \in U$
such that $f(x_U)\in [1/2, 1]$ and the points $f(x_U)$ are all
distinct.
Put $A=\{ x_U : U \in {\mathcal U} \}$, then $p \in \overline A$ and
so there exists a sequence $\{x_n : n<\omega \} \subset A$
converging to $p$, contradicting Lemma \ref{1.3'}.
\end{proof}

There is a variant of this result in which the assumption that
$Y$ be Tychonoff is relaxed to $T_3$, however the assumption on
the character of the point $p$ is more stringent.
Its proof will make
use of the following (probably well-known) lemma:

\begin{lemma}
\label{1.14} Let $Z$ be an infinite connected
regular space, then any non-empty open subset $G$ of $Z$
is uncountable.
\end{lemma}

\begin{proof} Choose a point $z \in G$ and an open
proper subset $V$ of $G$ with $z\in V \subset \overline V \subset G$.
If $G$ would be countable then, as a countable
regular space, $G$ would be Tychonoff, and so there would be a continuous
function $f: G \to [0,1]$ such that $f(z)=1$ and
$f$ is identically zero on $G-V$. Extend $f$ to a function
$\overline f : Z \to [0,1]$ by putting $\overline f(y)=0$ if $y\in
Z-G$. Then $\overline f$ is continuous and hence $\overline {f}(Z)$ is
also connected.
Consequently we have $f(G) = \overline{f}(Z)=[0,1]$ implying
that $|G|\ge |[0,1]|>
\omega$, and so contradicting that $G$ is countable.\end{proof}

\begin{theorem}
\label{1.15} Let $f$ be a preserving
function from a topological space $X$ into a regular space $Y$ and
let $p \in X$ be a Fr\` echet point of local connectivity with
character $\le \omega_1$. Then $f$ is continuous at $p$.
\end{theorem}

\begin{proof} Assume $f$ is discontinuous at the
point $p \in X$. As $p$ is a a Fr\` echet point, there is a
sequence $x_n \to p$ such that $f(x_n)$ does not converge to
$f(p)$. Taking a subsequence if necessary, we can suppose by Lemma
\ref{1.3'} that $f(x_n)=y \ne f(p)$ for all $n<\omega$.

Choose now an open neighbourhood $V$ of the point $y\in Y$ with
$f(p) \not \in \overline V$. Let $\mathcal U$ be a neighbourhood base
of the point $p$ in $X$ such that $|{\mathcal U}|\le \omega_1$ and the
elements of $\mathcal U$ are connected.

Choose now points $x_U$ from the sets $U \in \mathcal U$ such that
$f(x_U)\in V$ and the points $f(x_U)$ are all distinct. This can
be accomplished by an easy transfinite recursion because for each
$U\in \mathcal U$ the set $f(U)$ is connected and infinite, hence
$f(U)\cap V$ is uncountable by the previous
lemma. Put $A=\{x_U : U \in {\mathcal U} \}$. Then $p \in \overline A$
and so there exists a sequence $\{y_n : n<\omega \} \subset A$
converging to $p$, contradicting Lemma \ref{1.3'}.
\end{proof}

We shall now consider some further topological properties and prove
several results about them saying that preserving functions are
sequentially continuous.
Since in a Fr\'echet point sequential continuity implies continuity, these
results are clearly relevant to McMillan's theorem. 
Their real significance, however, will only become clear in the following
two sections.

\begin{definition}
\label{1.16} A point $x$ in a
topological space $X$ is called a {\it sequentially connectible
(in short: SC)} point, if $x_n \in X$, $x_n \to x$ implies that there are an
infinite subsequence $\langle x_{n_k} : k<\omega \rangle $ and a
sequence $\langle C_k :k<\omega \rangle $ consisting of connected
subsets of $X$ such that $\{x_{n_k},x\} \subset C_k$ for all $k<\omega$,
(i.e.\ $C_k$ ``connects'' $x_{n_k}$ with $x$, this explains the
terminology,)
moreover $C_k \to x$, i.e.\ every neighbourhood of
the point $x$ contains all but finitely many $C_k$'s.
A space $X$ is called
an $SC$ space if all its points are $SC$ points.
\end{definition}

\begin{remark}
\label{1.17}
It is clear that the SC property is closely related to local connectivity.
Let us say that a point $x$ in
space $X$ is a strong local connectivity point if it
has a neighbourhood base $\mathcal B$ such that the
intersection of an arbitrary (non-empty) subfamily of $\mathcal B$ is
connected. For example, local connectivity points
of countable character
or any point of a connected linearly ordered space has this
property.

We claim that {\it if $x$ is a strong local connectivity point of $X$
then $x$ is an $SC$ point in $X$.} Indeed,
assume that $x_n \to x$ and for every $n
\in \omega$ let $C_n$ denote the intersection of all those members of
$\mathcal B$ which contain both points $x_n$ and $x$.
(As the sequence $\{ x_n : n<\omega \}$ converges to $x$, we can
suppose that some element $B_0 \in {\mathcal B }$ contains all the
$x_n$'s.) Then $\{x,x_n\} \subset C_n$, moreover $C_n \to x$. Indeed,
the latter holds because if
$x\in B \in {\mathcal B}$ then, by definition, $x_n \in B$ implies
$C_n \subset B$. \easytosee
\end{remark}

The $SC$ property does not imply local connectivity. (If every convergent
sequence is eventually constant then the space is trivially $SC$.) 
However, the following simple lemma shows that if there are
``many'' convergent sequences then such an implication is valid.

\begin{lemma}
\label{1.18} 
Let $x$ be a both Fr\`echet and $SC$ point in a space $X$. 
Then $x$ is also a point of local connectivity in $X$.
\end{lemma}

\begin{proof}
Let $G$ be any open set containing $x$ and set $$H = \cup \{\,C
:\, x \in C \subset G\,\,\, \hbox{and}\,\, C\,\, \hbox{is
connected} \,\}.$$ We claim that (the obviously connected) set $H$
is a neighbourhood of $x$. Indeed, otherwise, as $x$ is a
Fr\`echet point, we could choose a sequence $x_n \to x$ from the
set $G-H$ while for every point $y \in G-H$ no connected set
containing both $x$ and $y$ is a subset of $G$, contradicting the
$SC$ property of $x$.
\end{proof}

If $SC$ holds globally, i.e.\ in an $SC$ space, then in the above result
the Fr\`echet property can be replaced with a weaker property that will
turn out to play a very important role in the sequel.

\begin{definition}
\label{1.19} 
A point $p$ in a topological space $X$ is called an {\it s} point if for
every family $\mathcal A$ of subsets of $X$ such that 
$p \in \overline {\bigcup \mathcal A}$ but $p \not \in \overline A$ for
all $A \in \mathcal A$ there is a sequence 
$\langle \langle x_n,A_n \rangle :n<\omega \rangle$ such that 
$x_n \in A_n \in \mathcal A$, the sets $A_n$ are pairwise distinct and
$\{x_n \}$ converges to some point $x \in X$ (that may be different from
$p$).
\end{definition}

A Fr\`echet point is evidently an $s$ point, moreover any point
that has a sequentially compact neighbourhood is also an $s$ point.
Other examples of $s$ points will be seen later.

\begin{theorem}
\label{1.20} 
Any $s$ point in a $T_3$ and $SC$ space is a point of local connectivity.
\end{theorem}

\begin{proof} 
Let $p$ be an $s$ point in the regular $SC$ space $X$ and let $G$ be an
open neighbourhood of $p$. 
We have to prove that the component $K_0$ of the point $p$ in $G$ is a
neighbourhood of $p$. 
Assume this is false and choose an open set $U$ such that $p \in U \subset
\overline U \subset G$.

Put
\begin{displaymath}
{\mathcal A }= \{K\cap \overline U : K \hbox{ is a component of
}G,\,\, K\not = K_0 \}.
\end{displaymath}

Then $p \in \overline {\bigcup \mathcal A}$ and $p \not \in \overline
A$ for $A\in \mathcal A$ (because a component of $G$ is relatively
closed in $G$), hence, by the definition of an $s$ point, there
exists a sequence $\{ \langle x_n, A_n \rangle : n<\omega \}$
such that $x_n \in A_n \in \mathcal A$, $x_n \to x$ for some $x\in X$
and if $A_n = K_n \cap \overline U$ then the components $K_n$ are
distinct. Note that $x\in \overline U \subset G$. As distinct
components are disjoint, we can assume that $x \not \in K_n$ for all
$n<\omega$. As $x$ is an $SC$ point, there are a connected
set $C$ and some $n<\omega$ such that $\{x, x_n \}\subset C\subset
G$. However, this is impossible, because then $K_n \cup C$ would
be a connected set in $G$ larger then the component
$K_n$.\end{proof}

The significance of the $SC$ property in our study of
continuity properties of preserving functions is revealed by the
following result.

\begin{theorem}
\label{1.21} 
A preserving function $f:X \to Y$ into a Hausdorff space $Y$ is
sequentially continuous at each $SC$ point of $X$.
\end{theorem}

\begin{proof} 
Let $x \in X$ be an $SC$ point and assume that $x_n \to x$ but $f(x_n)$
does not converge to $f(x)$ for a sequence $\{x_n : n< \omega \}$ in $X$.
We can assume by Lemma \ref{1.3'} that $f(x_n)=y\not = f(x)$ for all
$n<\omega$.
Choose an open neighbourhood $V$ of $y$ in $Y$ such that 
$f(x)\not \in \overline V$.

As $x$ is an $SC$ point in $X$, we can also assume that there is a
sequence of connected sets $C_n$ such that $C_n \to x$ and $x,x_n \in C_n$
for $n<\omega$. 
We can now define a sequence $z_n \in C_n$ such that $f(z_n) \in V$ and
the points $f(z_n)$ are all distinct. 
Indeed, assume $n<\omega$ and the points $z_i$ are already defined for
$i<n$ in this way.
As $f(C_n)$ is connected and $f(C_n)\cap V$ is a non-empty open proper
subset, this intersection $f(C_n)\cap V$ is not closed and hence is
infinite. 
Consequently there exists a point $z_n \in C_n$ with 
$f(z_n) \in f(C_n)\cap V - \{f(z_i) :i<n \}$. 
But then the sequence $\{z_n \}$ contradicts Lemma \ref{1.3'}. 
\end{proof}

\begin{corollary}
\label{1.22.} 
Let $f$ be a preserving function from a topological space $X$ into a
Hausdorff space $Y$ and let $p$ be both an $SC$ point and a 
Fr\`echet-point in $X$. 
Then $f$ is continuous at $p$. \easytosee
\end{corollary}

The following example (due to E. R. McMillan \cite{McM}) yields a locally 
connected $SC$ space with a discontinuous preserving 
function. (Compare this with Theorem \ref{1.21}.)

\begin{example}
\label{1.23} Take $\omega_1$ copies of the interval $[0,1]$ and
identify the $0$ points. We get in this way a ``hedgehog''
$X=\{0\}\cup \{R_{\xi}:\xi \in \omega_1 \}$, where the spikes
$R_{\xi}=(0,1] \times \{\xi\}$ are disjoint copies of the half
closed interval $(0,1]$. A basic neighbourhood of a point $x\in
R_{\xi}$ is an open interval around $x$ in $R_{\xi}$. A basic
neighbourhood of $0$ is a set of the form $\{0\} \cup \bigcup
\{J_{\xi}:\xi <\omega_1 \}$, where each $J_{\xi}$ is a non-empty
initial interval of $R_{\xi}$ and $J_{\xi}=R_{\xi}$ holds for all
but countably many ordinals $\xi$.

It is easy to see that in this way we get a locally connected
Tychonoff $SC$ space $X$. The function $f: X \to\,[0,1] $ defined by
$f((x,\xi))=x, f(0)=0$ is preserving but not continuous at the
point $0$ because every neighbourhood of $0$ is mapped onto
$[0,1]$. \easytosee
\end{example}

The next example is locally connected, hereditary Lindelof,
$T_{3{\frac 1 2}}$, of countable tightness and with a preserving
function given on it that is not even sequentially continuous.
It follows that this space is not an $SC$ space.

\begin{example}
\label{1.24} The underlying set of our space $X$
consists of a point $p$, of a sequence of points $p_n$ for
$n<\omega$ and countably many arcs $\{I(n,m) : m< \omega\}$
with disjoint
interiors connecting the points $p_n$ and $p_{n+1}$ for every
$n<\omega$.

If $x$ is an inner point of some arc, then its basic
neighbourhoods are the open intervals around it on the arc. The
basic neighbourhoods of a point $p_n$ are the unions of initial
(or final) segments of the arcs containing $p_n$. Finally, basic
neighbourhoods of $p$ are the sets which contain all but finitely
many $p_n$'s together with their basic neighbourhoods and for any
two consecutive $p_n$'s in the set all but finitely many of the
arcs $I(n,m)$. Note that the subspace $X-\{p\}$ can be realised
as a subspace of the plane, hence it is easy to check that $X$
has the above stated properties.

Now let $f:X\to [0,1]$ be defined as follows: $f(p)=0$, $f(p_n)=0$ if
$n$ is even, $f(p_n)=1$ if $n$ is odd, and $f$ is continuous on
each arc $I(n,m)$.
Then $f$ is not sequentially continuous at $p$ as shown by
the sequence$\{p_n :
n \hbox{ odd }\}$ converging to p, but it is preserving. Indeed, an
infinite sequence whose members are from the interiors of
different arcs is closed discrete and so a compact subset of $X$
can meet only finitely many open arcs. It follows then that its $f$-image
is the union of finitely many compact subsets of [0,1].
Moreover, if a connected
set contains both $p$ and some other point, then it also contains an
arc $I(n,m)$, and thus its image is the whole $[0,1]$. \easytosee
\end{example}

In the rest of this section we shall consider a slight weakening
of the sequential continuity property that comes up naturally in
the proof of McMillan's theorem or Theorems \ref{2.3} and \ref{2.4} below.

\begin{definition}
\label{1.25} A function $f:X\to Y$ is
said to be {\it weakly sequentially continuous} at the point $x$
if $f(x_n) \to f(x)$
whenever $x_n \to x$ in $X$ and
$f$ is not locally constant at $x_n$ for all $n< \omega$.
\end{definition}

We shall give below two types of points in which preserving functions
turn out to be weakly sequentially continuous.
In the next section then these will be used to
yield ``real'' continuity of
preserving functions on some interesting classes of spaces.

\begin{definition}
\label{1.26} A point $x$ in a
space $X$ is called an {\it
inflatable point} if $x_n \to x$ with $x_n \not = x$ for all $n<\omega$
implies that there are a subsequence $\{x_{n_k} : k <
\omega \}$ and neighbourhoods $U_k$ of $x_{n_k}$ for $k<\omega$
such that $U_k \to x$ (i.e.\ every neighbourhood of $x$ contains
all but finitely many $U_k$'s). The space $X$ will be called
inflatable if all its points are.
\end{definition}

It is obvious that any generalized ordered space is inflatable.

\begin{theorem}
\label{1.27} 
Any preserving function $f:X\to Y$ from a locally connected space $X$ into
a $T_2$-space $Y$ is weakly sequentially continuous at an inflatable point
$x$.
\end{theorem}

\begin{proof} 
Assume, indirectly, that $x_n \to x$ but $f(x_n)$ does not converge to
$f(x)$, while $f$ is not locally constant at $x_n$ for all $n<\omega$. 
By Lemma \ref{1.3'} we can assume that $f(x_n)=y \not = f(x)$ for all
$n<\omega$. 
Choose an open neighbourhood $V\subset Y$ of $y$ such that 
$f(x) \not \in \overline V$. 
As $x$ is inflatable, we may also assume to have open sets $U_n$ with 
$x_n \in U_n$ such that $U_n \to x$. 
Using Lemma \ref{1.7} we can recursively choose points $z_n \in U_n$ with
distinct $f$-images such that $f(z_n) \in V$ for all $n<\omega$,
contradicting Lemma \ref{1.3'} again.
\end{proof}

The other property is both a weakening of the Fr\`echet property and a
variation on the $s$ property.

\begin{definition}
\label{1.28} 
We call a point $x$ in a space $X$ a {\it set-Fr\'echet} point if whenever
$A = \bigcup\{A_n : n<\omega \}$ with $x \in \overline A$ but 
$x \not \in \overline {A_n}$ for all $n < \omega$ then there is a sequence
$\{x_n \} \subset A$ such that $x_n \to x$. 
Of course, a space is set-Fr\'echet if all its points are.
\end{definition}

\begin{theorem}
\label{1.29} 
Let $f$ be a preserving function from a locally connected $T_2$-space $X$
into the interval $[0,1]$. 
Then $f$ is weakly sequentially continuous at every set-Fr\'echet point
$x$ of $X$.
\end{theorem}

\begin{proof} 
Assume $x_n \to x$ but $f(x_n)$ does not converge to $f(x)$, moreover $f$
is not locally constant at each $x_n$ for $n<\omega$. 
By Lemma \ref{1.3'}, we can assume that $f(x)=1$ and $f(x_n)=0$ for all
$n<\omega$.
Note that then for any connected neighbourhood G of any point $x_n$ the
image $f(G)$ is a proper interval containing $0$.
For every $n<\omega$ choose an open sets $U_n$ such that $x_n \in U_n$ and
$x \not \in \overline U_n$ and put $A_n = \{z \in U_n:0<f(z)<1/n+1 \}$. 
By Lemma \ref{1.7} the conditions in the definition of a set-Fr\'echet
point are satisfied for the sets $A_n$ so there is a sequence of points
$z_n \in A= \bigcup\{A_n : n<\omega \}$ converging to $x$. 
It is immediate that $f(z_n)$ converges to $0\not = 1 = f(x)$ while the
set $\{f(z_n):n<\omega \}$ is infinite because $f(z_n) \not = 0$ for all
$n$, contradicting Lemma \ref{1.3'}.
\end{proof}

\section{From sequential continuity to continuity}

The aim of this section is to prove that if a locally connected
space $X$ fulfills one of the ``convergence-type'' conditions of the
first section (i.e.\ $X$ is an $SC$-space or it is inflatable or
set-Fr\'echet) and $f:X \to Y$ is a preserving function then,
assuming in addition
appropriate separation axioms for $X$ and $Y$, $f$ is
continuous at every $s$-point of $X$. The proofs of these theorems,
just like their formulations,
are very similar.

\begin{theorem}
\label{2.1} A preserving function $f:
X \to Y$ from a locally connected $SC$-space $X$ into a regular
space $Y$ is continuous at every $s$-point of $X$.
\end{theorem}

\begin{proof}
Assume indirectly that $f$ is not continuous
at the $s$-point $p\in X$. Then there exists a closed set $F \subset Y$ such
that $p\in \overline {f^{-1}(F)}$ but $f(p) \not \in F$. Choose an
open set $V\subset Y$ such that $F \subset V$ and $f(p) \not \in \overline V$.

Let $\mathcal K$ be the family of the components of $f^{-1}(\overline V)$ and put
${\mathcal A} = \{K\cap f^{-1}(F) : K \in {\mathcal K} \}$.
As $p \not \in \overline K$
for $K \in \mathcal K$ by Lemma \ref{1.4} and $p\in \overline {f^{-1}(F)}$,
the conditions
given in the definition of an $s$ point are fulfilled for the family
${\mathcal A}$.
Thus there is a sequence $\{x_n : n<\omega\} \subset f^{-1}(F)$ such that
$x_n \to x$ for some $x\in X$ and if $K_n$ is the component of $x_n$ in
$f^{-1}(\overline V)$, then $K_m \not = K_n$ for $m \not = n$.

As the components $K_n$ are pairwise disjoint, we can suppose
that $x \not \in K_n$ for all $n<\omega$. It follows that if $C$ is
a connected set which contains both $x$ and some $x_n$ then $C\not
\subset f^{-1}(\overline V)$, because otherwise $K_n \cup C$
would be a connected subset of $f^{-1}(\overline V)$ strictly
larger than the component $K_n$, a contradiction. Hence, using that $x$ is an
$SC$-point, we may assume to have
a sequence $C_n$ of connected sets such that
$C_n \to x$ and $C_n \not \subset f^{-1}(\overline V)$.
We can choose points $z_n \in C_n - f^{-1}(\overline V)$
for all $n < \omega$, then $z_n \to x$ and $f(z_n) \not \in
\overline V$. But by
Theorem \ref{1.21} $f$ is sequentially continuous, consequently
$f(x) = \lim f(z_n) \in Y-V$. On the other hand,
$f(x_n)\in F$ for $n<\omega$ and so, using again the sequential
continuity of $f$ at the point $x$, we get that $f(x)=\lim
f(x_n)\in F$, a contradiction.\end{proof}

The proofs of the other two analogous theorems for inflatable and 
set-Fr\'echet spaces, respectively, make essential use of the following
lemma:

\begin{lemma}
\label{2.2}Assume that $f:X \to Y$ is
a preserving and weakly
sequentially continuous function from a locally
connected $T_3$-space $X$ into a $T_3$-space $Y$ and $f$ is not continuous
at some $s$-point of $X$. Then there are two sets $F\subset V \subset
Y$, $F$ closed and $V$ open in $Y$ and a convergent sequence $x_n \to x$
in $X$ such that for all $n<\omega$ we have $x_n \not = x,\,
f(x_n)\in F$ but $f(U)\not \subset \overline V$
whenever $U$ is a neighbourhood of $x_n$.
It follows that $f$ is not locally constant at the points $x_n$ and $x$.
\end{lemma}

\begin{proof} Assume $f$ is not continuous at
the $s$-point $p\in X$, then there exists a closed set $F \subset
Y$ such that $p\in \overline {f^{-1}(F)}$ but $f(p) \not \in F$.
Choose an open set $V\subset Y$ such that $F \subset V$ and $f(p)
\not \in \overline V$.

Put $$B=\{x\in f^{-1}(F): f \hbox{ is not locally constant at } x\},$$
then $p \in \overline B$ by Lemma \ref{1.6}. Now let 
$$H=\{x \in f^{-1}(F): f(U) \not \subset \overline V \hbox{ for every
neighbourhood } U \hbox{ of } x \},$$ 
clearly $H \subset B$. 
We assert that $p \in \overline H$ as well. 
To show this, fix a closed neighbourhood $W$ of $p$. 
Let $\mathcal K$ be the family of the components of $f^{-1} (\overline V)$
and put ${\mathcal A} = \{K\cap B\cap W : K \in {\mathcal K} \}$. Since $p
\not \in \overline K$ for $K \in \mathcal K$ by Lemma 1.4, the conditions
in the definition of an $s$-point are satisfied for $p$ and ${\mathcal
A}$. Thus there is a sequence $\{y_n : n<\omega\} \subset B\cap W$ such
that $y_n \to y$ for some $y\in X$ and if $K_n$ is the component of $y_n$
in $f^{-1}(\overline V)$, then $K_m \not = K_n$ if $m \not = n$.

We claim that $y\in W\cap H$. As $W$ is closed, trivially $y\in W$. The
sets $K_n$ are disjoint so we can assume that $y \not \in K_n$ for all
$n<\omega$. Using that $f$ is weakly sequentially continuous and $y_n \in
B$, we get that $f(y)=\lim f(y_n) \in F$, hence $y \in B$. We have yet to
show that $f(U)\not \subset \overline V$ if $U$ is any neighbourhood of
$y$. By local connectivity, we can suppose that $U$ is connected. But the
connected set $U$ meets (infinitely many) distinct components of
$f^{-1}(\overline V)$, so indeed $U \not \subset f^{-1}(\overline V)$.

Now applying the $s$-point property of $p$ to the family ${\mathcal A} =
\{\{x\}: x \in H \}$ there is an infinite sequence of points $x_n \in H$
converging to some point $x$, completing the proof.
\end{proof}

\begin{theorem}
\label{2.3} 
A preserving function from a locally connected and inflatable $T_3$ space
$X$ into a $T_3$ space $Y$ is continuous at every $s$-point of $X$.
\end{theorem}

\begin{proof} 
Assume, indirectly, that the preserving function $f:X \to Y$
is not continuous at an $s$-point of $X$. By Theorem \ref{1.27} $f$ is
weakly sequentially continuous, hence we can apply the preceding lemma and
choose appropriate sets $F$, $V$ in $Y$ and points $x_n$ and $x$ in $X$.
As $X$ is inflatable and locally connected, we can also assume that there
is a sequence of connected open sets $U_n$ such that $U_n \to x$ and $x_n
\in U_n$ for $n<\omega$. Then for all $n$ we have $f(U_n) - \overline V
\not = \emptyset$, hence by Lemma \ref{1.14} these sets are uncountable.
Consequently we can recursively select another sequence of points $y_n \in
U_n$ (and so converging to $x$) such that $f(y_n) \in Y-\overline V$ and
the $f(y_n)$'s are pairwise distinct. But then the sequence $f(y_n)$ does
not converge to $f(x)$ because $f(x)=\lim f(x_n)\in F \subset V$ by the
weak sequential continuity of $f$, contradicting Lemma \ref{1.3'} again.
\end{proof}

\begin{corollary}
\label{2.4} 
If $X$ is locally compact, locally connected, and monotonically normal
(in particular if $X$ is a locally connected LOTS) then $Pr(X,T_3)$ holds.
\end{corollary}

\begin{proof} 
See \cite[theorem 3.12]{Juhasz} for a proof that a (locally) compact,
monotonically normal space is both inflatable and radial. Consequently, it
is also locally sequentially compact, and thus an $s$-space.
\end{proof}

\begin{theorem}
\label{2.5} 
A preserving function from a locally connected and set-Fr\`echet $T_3$
space $X$ into a $T_3{\frac 1 2}$ space is continuous at every $s$-point
of $X$.
\end{theorem}

\begin{proof}
By Lemma \ref{1.2}, it suffices to prove this for preserving functions
$f:X \to [0,1]$. Assume, indirectly, that the function $f$
is not continuous at some $s$-point $p\in X$.
Then, by Theorem \ref{1.29}, $f$ is weakly sequentially
continuous, hence we may again apply Lemma \ref{2.2} to choose
appropriate sets $F$ and $V$ in $Y = [0,1]$
and points $x_n$ and $x$ in $X$. There is a
continuous function $g:[0,1] \to [0,1]$ which is identically $1$
on $F$ and $0$ on $[0,1]-V$. Then the composite function $h=gf:X
\to [0,1]$ is also preserving, $h(x_n)=1$ for all $n$, and the $h$-image of
any neighbourhood of a point $x_n$ is the whole interval
$[0,1]$.

Let us choose open sets $G_n$ for $n<\omega$ such that $x_n \in G_n$ and
$x\not \in \overline G_n$. If $A_n = G_n \cap
f^{-1}((0,{\frac 1 n}))$ and $A=\bigcup A_n$, then $x \not \in
\overline A_n$ but $x \in \overline A$. So, as $X$ is a
set-Fr\`echet space, there is a sequence of points $y_n \in A$ converging to
$x$. But then the set $\{h(y_n): n<\omega \}$ is infinite and
$h(y_n) \to 0 \not = 1 =f(x)$, contradicting Lemma \ref{1.3'}.
\end{proof}

\begin{corollary}
\label{2.6} 
If $X$ is a locally connected, locally countably compact, and
set-Fr\`echet $T_3$-space then $Pr(X,T_3{\frac 1 2})$ holds.
\end{corollary}

\begin{proof} 
It is enough to note that a countably compact set-Fr\`echet space is also
sequentially compact and so an $s$-space.
\end{proof}

\section{ Some product theorems}

The aim of this section is to prove that an arbitrary product $X$ of
certain ``good'' spaces has the property $Pr(X,T_3)$. To achieve this, we
shall make use of Theorem \ref{2.1}. It is well-known that the product of
connected and locally connected spaces is locally connected, moreover a
very similar argument (based on the productivity of connectedness) implies
that the product of connected $SC$ spaces is $SC$. Hence two of the
assumptions of Theorem \ref{2.1} are productive if the factors are
connected. Nothing like this can be expected, however, about the third
assumption of Theorem \ref{2.1}, namely the $s$-property. To make up for
this, we are going to consider a stronger property that is at least
countably productive, and use this stronger property to establish what we
want, first for $\Sigma$-products and then for arbitrary products. Now,
this stronger property will require the existence of a winning strategy
for player {\bf I} in the following game.

\begin{definition}
\label{D3.1}
Fix a space $X$ and a point $p\in X$. The game $G(X,p)$
is played by two players {\bf I} and {\bf II} in $\omega$ rounds.
In the n-th round first {\bf I} chooses a neighbourhood $U_n$ of $p$ and then
{\bf II} chooses a point $x_n \in U_n$. {\bf I} wins if the produced
sequence $\{x_n : n<\omega \}$ has a convergent subsequence,
otherwise {\bf II} wins.

We shall say that $X$ is {\it winnable} at the point $p$ if {\bf
I} has a winning strategy in the game $ G(X,p) $.
$X$ is {\it winnable} if $ G(X,p)$ is winnable for
all $p\in X$. Note that, formally, a winning strategy for {\bf I} is a map
$ \sigma : X^{<\omega} \to {\mathcal V}(p)$, where ${\mathcal V}(p)$ is
the family of neighbourhoods of $p$, such that if $ \langle x_n : n<\omega
\rangle $
is a sequence played following $\sigma$, i.e.\ $x_n \in
\sigma \langle x_0,x_1,...,x_{n-1}\rangle$ for all $n<\omega$, then
$ \langle x_n : n<\omega \rangle $ has a convergent subsequence.
\end{definition}

\begin{lemma}
\label{3.1} {\it A winnable point $p$ of a
space $X$ is always an $s$ point.}
\end{lemma}

\begin{proof} 
Let $\sigma$ be a winning strategy of {\bf I} in the game
$G(X,p)$ and fix a family $\mathcal A$ of subsets of $X$ with $p \in
\overline {\bigcup \mathcal A}$ but $p \not \in \overline A$ for all $A
\in \mathcal A$. Let us play the game $G(X,p)$ in such a way that {\bf I}
follows $\sigma$ and assume that the first $n$ rounds of the game have
been played with the points $x_i\in U_i$ and the distinct sets $A_i \in
\mathcal A$ with $x_i \in A_i$ chosen by player {\bf II} for $i<n$. Let
$U_n$ = $\sigma \langle x_0,x_1,...,x_{n-1}\rangle$ be the next winning
move of {\bf I}, then {\bf II} can choose a set $A_n \in \mathcal A$ with
$A_n \cap (U_n - \cup_{i<n} \overline A_i )\not = \emptyset$ and then pick
$x_n \in A_n \cap (U_n - \cup_{i<n} \overline A_i )$ as his next move. But
player {\bf I} wins hence, a suitable subsequence of $\{x_n : n<\omega
\}$, whose members were picked from distinct elements of $\mathcal A$,
will converge.
\end{proof}

In order to prove the desired product theorem for winnable spaces we shall
consider {\it monotone strategies} for player {\bf I}. 
A strategy $\sigma$ of player {\bf I} is said to be {\it monotone} if for
every subsequence $\langle x_{i_0},...,x_{i_{r-1}}\rangle$ of a sequence
$\langle x_0,x_1,...,x_{n-1} \rangle $ of points in $X$ we have
$$\sigma\langle x_0,x_1,...,x_{n-1} \rangle \subset \sigma\langle
x_{i_0},...,x_{i_{r-1}}\rangle .$$

\begin{lemma}
\label{3.2} 
If player {\bf I} has a winning strategy in the game 
$G \langle X,p \rangle$ then he also has a monotone winning strategy.
\end{lemma}

\begin{proof}
Let $\sigma$ be a winning strategy for player {\bf I}; we define
a new strategy $\sigma_0$ as follows : for any
sequence $s= \langle x_0, x_1,...,x_{n-1} \rangle $ put
$$\sigma_0(s)=\bigcap\{\sigma \langle x_{i_0},...,x_{i_{r-1}}\rangle :
0\le i_0<i_1<...<i_{r-1}<n \}.$$

The function $\sigma_0$ is clearly a monotone winning strategy for {\bf I}.
\end{proof}

It is easy to see that if {\bf I} plays using
the monotone strategy $\sigma_0$ then any infinite
subsequence of the sequence chosen by player {\bf II} is also a
win for {\bf I}, i. e. has a convergent subsequence.
Another important property of a monotone winning strategy $\sigma_0$ is the
following: If $\langle x_n : n < \omega \rangle$ is
a sequence of points in $X$ such that
we have
$x_n \in \sigma_0 \langle x_0,x_1,...,x_{n-1}\rangle$
only for $n \ge m$ for some fixed $m < \omega$ then this is still a winning
sequence for {\bf I}. Indeed, this holds because for every $n \ge m$
we have $$x_n \in \sigma_0 \langle x_0,x_1,...,x_{n-1}\rangle
\subset \sigma_0 \langle x_m,x_{m+1},...,x_{n-1} \rangle$$
by monotonicity, hence the ``tail'' sequence 
$\langle x_n : m \le n < \omega \rangle$ is produced by a play of the game
where {\bf I} follows the strategy $\sigma_0$.

For a family of spaces $\{X_s : s \in S \}$ and a fixed point $p$
(called the {\it base point}) of the product $X = \prod \{X_s : s
\in S \}$ let $T(x)$ denote the support of the point $x$ in $X$:
this is the set $\{s \in S : x(s) \not = p(s) \}$. Then
$\Sigma(p)$ (or simply $\Sigma$ if this does not lead to
misunderstanding) denotes the $\Sigma$-product with base point
$p$: it is the subspace of $X$ of the points with countable
support, i.e.
$$\Sigma = \{ x \in X: |T(x)|\le \omega \}.$$

In the proof of the next result we shall use two lemmas. The
first one is an easy combinatorial fact:

\begin{lemma}
\label{3.3} {\it Let $
\langle H_k : k<\omega \rangle $ be a sequence of countable
sets, then for every $n<\omega$ there is a finite set
$F_n$ depending only on the first $n$ many sets
$\langle H_k : k<n \rangle $ such that
$F_n \subset \bigcup_{i<n}H_i$, $F_n \subset F_{n+1}$ and $\bigcup
F_n = \bigcup H_n$.}
\end{lemma}

\begin{proof} 
Fix an enumeration $H_i= \{x(i,j): j<\omega \}$ of the set $H_i$ for all
$i<\omega$ and then let $F_n=\{x(i,j):i,j < n \}$.
\end{proof}

The second lemma is about the sequences which play a crucial role
in the games $G\langle X,p\rangle$. Let us call a sequence $\{x_n\}$ in the
space $X$ {\it good} if every infinite subsequence of it has a
convergent subsequence.

\begin{lemma}
\label{3.4} If $x_n \in X=\prod \{X_i :
i<\omega \}$ for $n<\omega$ and $\{x_n(i):n<\omega \}$ is a good
sequence in $X_i$ for all $i\in \omega$ then $\{x_n \}$ is a good
sequence in $X$.
\end{lemma}

\begin{proof} We shall prove that if $N$ is any infinite
subset of $\omega$ then there is in $X$ a convergent
subsequence of $\{x_n : n\in N \}$.

We can choose by recursion on $k< \omega$ infinite sets $N_k$ such that
$N_{k+1}\subset N_k\subset N$ and $\{x_n(k):n\in N_k \}$ converges
to a point $x(k)$ in $X_k$. Then there is a
diagonal sequence $\{n_k:k<\omega
\}$ such that $n_k \in N_k$ and $n_k < n_{k+1}$ for all $k<\omega$.
The sequence $\{n_k:k<\omega \}$ is eventually contained in $N_i$
, hence $x_{n_k}(i)\to x(i)$ in $X_i$,
for all $i<\omega$. It follows that
$\{x_{n_k}\}$ is a convergent subsequence of $\{x_n :
n\in N \}$ in $X$.\end{proof}

The following result says a little more than that winnability is a countably
productive property.

\begin{lemma}
\label{3.5} 
Let $p \in X=\prod \{X_s : s \in S \}$ and suppose that 
$G\langle X_s,p(s) \rangle $ is winnable for every $s\in S$. 
Then $ G\langle \Sigma(p),p \rangle $ is also winnable.
\end{lemma}

\begin{proof} We have to construct a winning strategy
$\sigma$ for player {\bf I} in the game $ G\langle \Sigma(p), p
\rangle $. By Lemma \ref{3.2}, we can fix a monotone winning
strategy $\sigma_s$ of {\bf I} in the game $
G\langle X_s, p(s) \rangle $ for each $s\in S$. Given a sequence
$\langle x_0,...,x_{n-1}\rangle \in [\Sigma (p)]^{<\omega},$
let $H_i$ denote the support of
$x_i$ and $F_n$ be the finite set assigned to the sequence $\{H_i\ : i < n\}$
as in Lemma \ref{3.3}. Now, if $\pi_s$ is the projection from the
product $X$ onto the factor $X_s$ for $s \in S$ then set

$$\sigma \langle x_i:i<n \rangle = \bigcap
\{\pi^{-1}_s (\sigma_s\langle x_i(s) :i<n\rangle): s\in F_n \rangle\}.$$

Now let $\langle x_n : n < \omega \rangle$ be a sequence of points
in $\Sigma(p)$ produced by a play of the game $ G\langle \Sigma(p), p
\rangle $ in which {\bf I} followed the strategy $\sigma$.
Then for every $s \in H = \bigcup H_n$ the sequence $\langle x_n (s):n
< \omega \rangle$ is a win for player {\bf I}
in the game $G\langle X_s, p(s) \rangle $ because there is an $m <
\omega$ with $s \in F_m$ and then $x_n(s) \in \sigma_s\langle x_i(s)
:i<n\rangle$ is valid for all $n \ge m$.
Consequently, by Lemma
\ref{3.4}, the sequence $\langle x_n |H : n<\omega \rangle$ has a
convergent subsequence in $\prod \{X_s : s \in H \}$, while
for $s \in S -H$ we have $x_n(s)=p(s)$ for all $n<\omega $, and so
$\langle x_n : n<\omega \rangle$ indeed has a
convergent subsequence in $\Sigma(p)$. \end{proof}

The following two statements both easily follow from the fact that any
product of connected spaces is connected.

\begin{lemma}
\label{3.6}
A $\Sigma$-product of connected $SC$ spaces is also an $SC$ space. 
\easytosee
\end{lemma}

\begin{lemma}
\label{3.7}
A $\Sigma$-product of connected and locally connected spaces is also
locally connected. 
\easytosee
\end{lemma}

We now have all the necessary ingredients needed to prove our main
product theorem.

\begin{theorem}
\label{3.8} Let $f:X =\prod \{X_s :
s\in S \} \to Y$ be a preserving function from a product of
connected and locally connected $SC$ spaces into a regular space $Y$.
If $p\in X$ and $ G\langle X_s,p(s) \rangle $ is winnable for all
$s\in S$ then $f$ is continuous at the point $p$ .
\end{theorem}

\begin{proof} Let $\Sigma$ denote the sigma-product
with base point $p$. Then, by Lemma \ref{3.5}, $ G\langle \Sigma,p
\rangle $ is winnable and so $p$ is an $s$ point in $\Sigma$.
Moreover, by Lemmas \ref{3.6} and \ref{3.7}, $\Sigma$ is also a locally
connected
$SC$ space. Hence Theorem \ref{2.1} implies that the restriction of the
function $f$ to the subspace $\Sigma$ of $X$ is continuous at $p$.

To prove that $f$ is also continuous at the point $p$ in $X$,
fix a neighbourhood $V$ of $f(p)$ in $Y$. As the restriction
$f|\Sigma$ is continuous at $p$, there is an elementary neighbourhood
$U$ of $p$ in the product
space $X$ such that $f(U\cap \Sigma) \subset V$. Since the
factors $X_s$ are connected and locally connected, we can assume
that $U$ and $U\cap \Sigma$ are also connected
and hence, by Lemma \ref{1.4}, we have
$$f(U)\subset f(\overline {U\cap\Sigma})\subset \overline {f(U\cap
\Sigma)} \subset \overline V.$$
The regularity of $Y$ then implies that $f$ is continuous at
$p$.\end{proof}

\begin{corollary}
\label{3.9} $Pr(X,T_3)$ holds whenever $X$
is any product of connected and locally connected, winnable $SC$
spaces. In particular, if $X=\prod\{X_s:s\in S \}$ where each
factor $X_s$ is either a connected linearly ordered space (with the
order topology) or a connected and locally connected first
countable space then $Pr(X,T_3)$ is valid.\qed
\end{corollary}

For the proof of the next Corollary we need a general fact about the
relations $Pr(X,T_i)$.

\begin{lemma}
\label{3.10} If $q \colon X \to Y$ is a
quotient mapping of $X$ onto $Y$ then, for any i, $Pr(X,T_i)$
implies $Pr(Y,T_i)$.
\end{lemma}

\begin{proof} Let $f \colon Y \to Z$ be a preserving
function into the $T_i$ space $Z$.
The function $fq \colon X \to Z$, as the composition of
a continuous (and so preserving) and of a preserving function is
also preserving, hence, by $Pr(X,T_i)$, it is continuous. But then
$f$ is continuous because $q$ is quotient.\end{proof}

\begin{corollary}
\label{3.11} Let $X=\prod\{X_s:s\in
S \}$ where all the spaces $X_s$ are compact, connected, locally
connected, and monotonically normal. Then $Pr(X,T_2)$ holds.
\end{corollary}

\begin{proof} 
It follows from the recent solution by Mary Ellen Rudin of Nikiel's
conjecture \cite{rudin}, combined with results of L.B.Treybig \cite{tr} or 
J.Nikiel \cite{nik}, that {\it every compact, connected, locally 
connected, monotonically normal space is the continuous image of a
compact, connected, linearly ordered space.} Hence our space $X$ is the
continuous image of a product of compact, connected, linearly
ordered spaces. But any $T_2$ continuous image of a compact $T_2$ space
is a quotient image (and $T_3$), hence Corollary \ref{3.9} and
Lemma \ref{3.10} imply our
claim. 
\end{proof}

Comparing this result with Corollary \ref{2.4}, the following question is
raised naturally.

\begin{problem}
\label{3.12}
Let $X$ be a product of locally compact,connected and locally connected
monotonically normal spaces. 
Is then $Pr(X,T_3)$ true?
\end{problem}

The following is result is mentioned here mainly as a curiosity.

\begin{corollary}
\label{3.13}
Let $X= \prod \{X_s : s\in S \} $ be a product of linearly ordered and/or
first countable spaces. 
Then the following are equivalent:
\begin{enumerate}[a)]
\item $Pr(X,T_3)$;
\item
$X$ is locally connected;
\item
the spaces $X_s$ are all locally connected and all but finitely
many of them are also connected.
\end{enumerate}
\end{corollary}

\begin{proof} a)$\Rightarrow$b) : Lemma \ref{1.1}.
\\
b)$\Rightarrow$c) : \cite[6.3.4]{Eng}
\\
c)$\Rightarrow$a) : Corollary \ref{3.9}.\end{proof}

\begin{remark}
\label{3.14} E. R. McMillan raised the
following question in \cite{McM}: does $Pr(X,T_2)$ imply that $X$ is a
$k$-space? We do not know the (probably negative) answer to this question,
however we do know that the answer is negative if $T_2$ is replaced by $T_3$
in it. Indeed, for instance
${\mathbb R}^{\omega_1}$ is not a $k$-space
(see e.g.\ \cite[exercise 3.3.E]{Eng},), but, by Corollary \ref{3.9},
$Pr({\mathbb R}^{\omega_1},T_3)$ is valid.
\end{remark}

\section{ The sequential and the compact cases}

The two examples given in section 1
of (locally connected) spaces on which there are
non-continuous preserving functions both lack the properties of
(i) sequentiality and (ii) compactness. Here (i) arises naturally
as a weakening of the Fr\`echet property from McMillan's theorem
while the significance of (ii) needs no explanation. This leads us
naturally to the following problem.

\begin{problem}
\label{4.1} Assume that the locally connected
space $X$ is (i) sequential and/or (ii) compact.
Is then $Pr(X,T_2)$ (or $Pr(X,T_{3{\frac 1 2}}))$ true?
\end{problem}

The answer in case (i) turns out to be
positive if we assume the $SC$ property instead of
local connectivity. The
following result reveals why local connectivity need not be assumed
in it.(Compare this also with Lemma \ref{1.18}.)

\begin{theorem}
\label{4.2} Any sequential $SC$ space
$X$ is locally connected.
\end{theorem}

\begin{proof} We have to prove that if $K$ is a
component of an open set $G\subset X$ then $K$ is open. Assume
not, then $X-K$ is not closed, hence as $X$ is sequential there
is a sequence $\{x_n \}\subset X-K$ such that $x_n \to x \in K$. Since
$X$ is an $SC$ space and $G$ is a neighbourhood of $x$
there is a connected set $C\subset G$ such
that $\{x, x_n\} \subset C$ for some $n<\omega$. But this is impossible
because then the connected set $K\cup C\subset G$ would be larger than the
component $K$ of $G$.\end{proof}

Now we give the above promised partial solution to Problem \ref{4.1} in 
case (i), i. e. for sequential spaces.

\begin{theorem}
\label{4.3} 
If $X$ is a sequential $SC$ space then $Pr(X,T_2)$ holds.
\end{theorem}

\begin{proof}
Let $f : X \to Y$ be a preserving map into a $T_2$ space $Y$. Since $X$ is
sequential it suffices to show that the function $f$ is sequentially
continuous but this immediate from Theorem \ref{1.21}.
\end{proof}

Our next result implies that a counterexample to Problem \ref{4.1} (in
either case) can not be very simple in the sense that discontinuity of a
preserving function cannot occur only at a single point (as it does in
both examples \ref{1.23} and \ref{1.24}). 
In order to prepare this result we first introduce a topological property
that generalizes both sequentiality and (even countable) compactness.

\begin{definition}
\label{4.4} 
$X$ is called a {\em countably $k$} space if for any set $A \subset X$
that is not closed in $X$ there is a countably compact subspace $C$ of $X$
such that $A \cap C$ is not closed in $C$.
\end{definition}

This condition means that the topology of $X$ is determined by its
countably compact subspaces. All countably compact and all $k$
(hence also all sequential) spaces are countably $k$.
It is easy to see that the countably $k$ property is always inherited
by closed subspaces and for regular spaces by open subspaces as well.

\begin{theorem}
\label{4.5} Let $X$ be countably
$k$ and locally connected and $Y$ be $T_3$, moreover let $f\colon X
\to Y$ be a preserving function. Then the set of points of
discontinuity of $f$ is not a singleton. So if $X$ is also $T_3$
then the discontinuity set of $f$ is dense in itself.
\end{theorem}

\begin{proof}
Assume, indirectly, that
$f$ is not continuous at $p \in X$
but it is continuous at all other points of $X$.
Then we can choose a closed set $F\subset Y$ with
$A=f^{-1}(F)$ not closed. Evidently, then $\overline A - A = \{p\}$.
As $A$ is not closed in $X$ and $X$ is countably $k$, there is
a countably compact set $C$ in $X$ such that $A\cap C$ is not
closed in $C$; clearly then $p \in C$
and $p$ is in the closure of $A\cap C$.

Let $V \subset Y$ be an open set with $F\subset V$ and $f(p) \not
\in \overline V$ and put $H=f^{-1}(V)$. Then $H$ is open in $X$
and contains the set $A$. For every component $L$ of $H$ we have
$p \not \in \overline L$ by
$f(p) \not \in \overline L \subset \overline V $ and
Lemma \ref{1.4}, hence $p \in \overline {A \cap C}$
implies that there are infinitely many components of $H$ that meet
$A \cap C$. Thus we may choose a sequence $\{L_n : n < \omega\}$
of distinct such components with points $x_n \in {L_n \cap A \cap C}$.

We claim that $x_n \to p$. As the $x_n$'s are chosen from the
countably compact set $C$, it is enough to prove
for this that if $x \ne p$ then $x$ is not an
accumulation point of the sequence $\{ x_n \}$.
If $x \not \in \overline A = A \cup \{p\}$ this is obvious so
we may assume that $x \in A \subset H$. But then, as $H$ is open and $X$ is
locally connected, the connected
component of $x$ in $H$ is a neighbourhood of $x$ that
contains at most one of the points
$x_n$.

The point $f(p)$ is not in the closure of the set $\{f(x_n) \colon
n<\omega \} \subset F$, hence, by Lemma \ref{1.3'}, we can suppose that
$f(x_n) = y \ne f(p)$ for each $n < \omega$.

Now we choose a sequence of neighbourhoods $V_n$ of $y$ in $Y$ with
$V_0 = V$ and $\overline{V_{n+1}} \subset V_n$ for all
$n<\omega$ and then put $U_n = f^{-1}(V_n)$.
Clearly $U_n$ is open in $X$ and $U_0 = H$, hence, as was noted above,
the closure of any connected set contained in $U_0$
contains at most one of the points $x_n$.

Next, let $K_n$ be the component of $x_n$ in $U_n$ for $n<\omega$.
We claim that every boundary point of $K_n$ is mapped by $f$ to
a boundary point of $V_n$, i. e. $$f(FrK_n) \subset FrV_n.$$
Indeed, $f(\overline{K_n}) \subset \overline{V_n}$
by Lemma \ref{1.4} (or by continuity at all points distinct from $p$).
Moreover, we have $$FrK_n = \overline{K_n} - K_n \subset FrU_n =
\overline{U_n} - U_n$$ because $K_n, U_n$ are open and $K_n$ is
a component of $U_n$,
and $f(FrU_n) \cap V_n = \emptyset$ because $U_n = f^{-1}(V_n)$.
Thus indeed $f(FrK_n) \subset \overline{V_n} - V_n = FrV_n$.

Every connected neighbourhood of $p$ meets
all but finitely many $K_n$'s hence also
$FrK_n$, consequently $K=\bigcup \{FrK_n : n<\omega
\}$ is not closed. So there exists a countably compact set $D$
with $D\cap K$ not closed in $D$. But the sets $FrK_n$ are
closed, thus $D$ must meet infinitely many of them, i. e. the set
$$N = \{ n < \omega : D \cap FrK_n \not = \emptyset \}$$
is infinite.
Let us choose a point $z_n$ from each nonempty $D \cap K_n$.

Again we claim that the
only accumulation point of the sequence $\{z_n : n \in N\}$ is $p$. Indeed,
if $x \not = p$ would be such an accumulation point, then $f(z_n)
\in \overline {V_n} \subset \overline {V_1}$ for all $n \in N - \{0\}$
would also also imply $f(x) \in
\overline {V_1} \subset V_0$. By continuity and local connectivity at $x$
then
there is a connected neighbourhood $W$ of $x$ with $f(W) \subset
{V_0}$. But then the set
$$W\cup \bigcup \{K_n : W \cap K_n \not =\emptyset \}$$
would be a connected subset of $U_0$ that has
$p$ in its closure, a contradiction.

Consequently, the sequence $\{z_n : n \in N\}$ must converge to $p$,
while $\{f(z_n) : n \in N\} \subset \overline V$ does not converge to $f(p)$,
contradicting Lemma \ref{1.3'} because $f(z_n) \in FrV_n $
for all $n \in N$ and the boundaries $FrV_n $ are pairwise disjoint, hence
the $f(z_n)$'s are pairwise distinct.

The last statement of the theorem now follows easily because
an isolated discontinuity of $f$ yields an open subspace of $X$
on which the restriction of $f$ has a single point of discontinuity,
although if $X$ is $T_3$ then
any open subspace of $X$ is both countably $k$ and locally connected.
\end{proof}

Noting that Problem \ref{4.1} really comprises three different questions,
and having shown above that, in a certain sense, it seems to be hard to
find counterexamples to any of these, we now turn our attention
to the case in which both (i) and (ii) are assumed. In this case
we can provide a
positive answer, at least consistently
and with the extra assumption that the cellularity
of the space in question is ``not too large''. In fact, what we
can prove is that if $2^{\omega} < 2^p$ then any locally connected,
compact $T_2$ space $X$ that is sequential
and does not contain a cellular family of size $p$ satisfies $Pr(X,T_2)$.
Of course, here $p$ stands for the well-known cardinal invariant of the
continuum whose definition is recalled below.

A set $H\subset \omega$ is called a
{\it pseudo intersection} of the family ${\mathcal A}\subset
[\omega]^{\omega}$ if $H$ is almost contained in every member of
${\mathcal A}$, i.e.\ $H-A$ is finite for each $A\in \mathcal A$. Then $p$
is the minimal cardinal $\kappa$ such that there exists a family
${\mathcal A}\subset[\omega]^{\omega}$
of size $\kappa$ which has the finite intersection property
but does not have an infinite pseudo intersection.
(Here the finite intersection property means that any finite
subfamily of ${\mathcal A}$ has {\it infinite} intersection.)

It is well-known (see e.g.\cite{GenTop}) that the cardinal $p$ is
regular, $\omega_1 \le p \le 2^{\omega}$ and
$2^{\kappa}=2^{\omega}$ for $\kappa < p$. The condition
``$2^p > 2^{\omega}$'' of our result is satisfied if $p=
2^{\omega}$ (hence Martin's axiom implies it), but it is also true if
$2^{\omega_1} > 2^{\omega}$.

Now, our promised consistency result on compact sequential spaces
will be a corollary of a ZFC result of somewhat technical nature.
Before formulating this, however, we shall prove two
lemmas that may have some independent interest in themselves.

\begin{lemma}
\label{4.6} Let $X$ be a compact $T_2$
space of countable tightness and $f : X \to [0,1]$ be a compactness
preserving map of $X$ into the unit interval. If $x \in X$ is a point
in $X$ and $[a,b]$ is a subinterval of $[0,1]$ such that for every
neighbourhood $U$ of $x$ we have $[a,b] \subset f(U)$ then for any
$G_{<p}$ set $H$ containing $x$ we also have
$[a,b] \subset f(H)$.
\end{lemma}

\begin{proof} Without loss of generality we may assume
that $H$ is closed. Now the proof will proceed by induction on $\kappa$
where $\omega \le \kappa < p$ and $H$ is a (closed) $G_{\kappa}$ set,
or equivalently, the character $\chi(H,X) = \kappa$. If $\kappa = \omega$ then
we can write $H=\bigcap \{G_n :n<\omega \}$ with $G_n$ open and
$\overline {G_{n+1} }\subset G_n$ for all $n<\omega$. Fix a countable
dense subset $\{c_n : n < \omega \}$ of $[a,b]$ and then pick $x_n \in G_n$
with $f(x_n) = c_n$, this is possible by assumption. Note that then every
accumulation point of the set $M = \{x_n : n < \omega\}$ is in $H$, hence
by Lemma \ref{1.3} we have $$[a,b] = f(M)' \subset f(M') \subset f(H).$$

Next, if $\omega < \kappa < p$ then we have
$x\in H = \bigcap \{S_{\xi} : \xi <
\kappa \}$, where $S_{\xi}\supset S_{\eta}$ if $\xi <\eta$ and the
$S_{\xi}$ are closed sets of character $<\kappa$. By induction, we have
$f(S_{\xi})\supset [a,b]$ for all $\xi < \kappa$, and we have to
prove that $f(H)\supset [a,b]$ as well. In fact, it suffices to show that
$f(H) \cap [a,b] \not = \emptyset$ because applying this to all
(non-singleton) subintervals
of $[a,b]$ we actually get that $f(H) \cap [a,b]$ is dense in $[a,b]$ while
$f(H)$ is also compact, hence closed.

We do this indirectly; assume $f(H) \cap [a,b] = \emptyset$ then we can
choose points $x_{\xi}\in S_{\xi} - H$ for all $\xi <\kappa$ such that the
images $f(x_{\xi})\in [a,b]$ are all distinct. Let $\bar x$ be a
complete accumulation point of the set $\{x_{\xi} : \xi <
\kappa \}$. Then $\bar x \in H$ and $t(X)=\omega$ implies that
there is a countable subset $A \subset \{x_{\xi} : \xi < \kappa
\}$ such that $\bar x \in \overline A - A$. Choose now a neighbourhood
base $\mathcal B$ of $H$ in $X$ of size $\kappa <p$. The family $\{A\cap B :
B \in {\mathcal B} \}\subset [A]^{\omega}$ has the finite intersection
property hence it has an infinite pseudo intersection $P \subset
A$, i. e. the set $P - B$ is finite for each $B \in \mathcal B$. This
implies that every accumulation point of $P$ is contained in $H$.
But $\overline P$ is compact, hence by Lemma \ref{1.3} we have
$$\emptyset \not = f(P)' \subset f(P') \cap [a,b] \subset f(H) \cap [a,b],$$
which is a contradiction.
\end{proof}

Before we give the other lemma, let us recall that for any space $X$
we use $\widehat c(X)$ to denote the smallest cardinal $\kappa$ such that
$X$ does {\it not} contain $\kappa$ disjoint open sets.

\begin{lemma}
\label{4.7} Let $f : X \to Y$ be a
connectivity preserving map from a locally connected space $X$ into
a $T_2$ space $Y$. Then for every $x \in X$ with $\chi(f(x),Y) < \widehat c(X)$
there is a $G_{<\widehat c(X)}$ set $H$ in $X$ such that $x \in H$ and if
$z \in H$ is any point of continuity of $f$ then $f(z) = f(x)$.
\end{lemma}

\begin{proof} Let $\kappa = \widehat c(X)$ and fix a
neighbourhood base $\mathcal V$ of the point $f(x)$ in $Y$ with
$|\mathcal V| < \kappa$. For every $V \in \mathcal V$ let us then set
$$G_V = \bigcup\{G : G\hbox{ is open in }X \,\hbox{and}\, f(G) \cap V =
\emptyset\}.$$
For every component $K$ of the open set $G_V$ we have $f(x) \not \in
\overline{f(K)}$ and therefore $x\not \in \overline{K}$ by Lemma \ref{1.4},
moreover the components of $G_V$ form a cellular family because $X$
is locally connected, hence their number is less than $\kappa$.
Consequently, $$H_V = \bigcap\{X - \overline K : K\, \hbox{is a component
of}\,\, G_V\}$$ is a $G_{< \kappa}$
set with $x \in H_V$ and $H_V \cap G_V = \emptyset$.

The cardinal $\kappa$ is regular (see e.g.\ \cite[4.1]{Juh}), hence
$H = \cap\{H_V : V \in \mathcal V\}$ is also a $G_{<\kappa}$ set that
contains the point $x$. Now, suppose that $z$ is a point of continuity
of $f$ with $f(z) \not = f(x)$. Then there is a basic neighbourhood
$V \in \mathcal V$ of $f(x)$ and a neighbourhood $W$ of $f(z)$ with
$V \cap W = \emptyset$, and there is an open neighbourhood $U$ of $x$ in
$X$ with $f(U) \subset W$. But then, by definition, we have $z \in U
\subset G_V$, hence $z \not \in H_V \supset H$. \end{proof}

\begin{theorem}
\label{4.8} Let $X$ be a locally connected
compact $T_2$ space of countable tightness. If, in addition, we also have
$|X| < 2^p$ and $\widehat c(X) \le p$ then $Pr(X,T_2)$ holds.
\end{theorem}

\begin{proof} 
Using Lemma \ref{1.2} it suffices to show that any
preserving function $f : X \to [0,1]$ is continuous. To this end, first
note that if $f$ is not continuous at a point $x \in X$ then the
oscillation of $f$ at $x$ is positive, hence , by local connectivity at
$x$ and because $f$ is preserving there are $0 \le a < b \le 1$ such that
$f(U) \supset [a,b]$ holds for every neighbourhood $U$ of $x$.
Consequently, by Lemma \ref{4.6} we also have $f(H) \supset [a,b]$
whenever $H$ is any $G_{<p}$ set containing the point $x$. In particular,
this implies that if the singleton $\{x\}$ is a $G_{<p}$ set
(equivalently, if the character of $x$ in $X$ is less than $p$) then $f$
is continuous at $x$.

On the other hand, by Lemma \ref{4.7}, for every point $x \in X$ there is
a closed $G_{<p}$ set $H_x$ with $x \in H_x$ such that for any point of
continuity $z \in H_x$ of $f$ we have $f(z) = f(x)$. We claim that $f$ is
constant on every such set $H_x$ and then, by the above, $f$ is continuous
at every point $x \in X$.

For this it suffices to show that $f$ has a point of continuity in every
(non-empty) closed $G_{<p}$ set $H$. Indeed, for any point $y \in H_x$
then the intersection $H_x \cap H_y$ contains a point of continuity $z$
for which $f(x) = f(z) = f(y)$ must hold. By the \v Cech-Pospi\v sil
theorem (see e.g.\ \cite[3.16]{Juh}) and by $|H| < 2^p$ there is a point 
$z \in H$ with $\chi(z,H) < p$ and so $\chi(z,X) < p$ as well, for $H$ is 
a $G_{<p}$ set in $X$. But we have seen above that then $z$ is indeed a
point of continuity of $f$.
\end{proof}

\begin{theorem}
\label{4.9}
Assume that $2^{\omega} < 2^p$ and X is a locally connected and sequential
compact $T_2$ space with $\widehat c(X) \le p$. Then $Pr(X,T_2)$ holds.
\end{theorem}

\begin{proof}
By a slight strengthening of some well-known results of Shapirovski (see
e.g.\ \cite[2.37 and 3.14]{Juh}), for any compact $T_2$ space $X$ we have
both $\pi\chi(X) \le t(X)$ and 
$$d(X) \le \pi\chi(X)^{<\widehat c(X)}.$$
Consequently, for our space $X$ we have 
$$d(X) \le \omega^{<p} = 2^{\omega}$$ 
and so by sequentiality $|X| \le 2^{\omega}$ as well. 
But this shows that all the conditions of Theorem \ref{4.8} are satisfied
by our space $X$. 
\end{proof}

To conclude, let us emphasize again that Lemma \ref{1.3}, i.e.\ the full
force of compactness preservation, as opposed to just the preservation of
the compactness of convergent sequences, was only used in this section
(cf.\ the remark made after \ref{1.3'}).

\section{The relation $Pr(X,T_1)$}

The main aim of this section is to prove that if $Pr(X,T_1)$ holds
and $X$ is $T_3$ then $X$ is discrete. Note the striking contrast
between $Pr(X,T_1)$ and $Pr(X,T_2)$: the latter holds for a large
class of (non-discrete) spaces (see Theorem \ref{1.8} or
Corollary \ref{3.10}).

Let us recall that the {\em cofinite
topology} on an underlying set $X$ is the coarsest $T_1$ topology on
$X$: the open sets are the empty set and the complements of the
finite subsets of $X$. It is not hard to see that such a space is
hereditarily compact and any infinite subset in it is connected.
Let us start with a result that gives several different characterizations
of $T_1$ spaces $X$ that satisfy $Pr(X,T_1)$.

\begin{theorem}
\label{5.1}
For a $T_1$ space $X$ the following conditions are
equivalent:
\begin{enumerate}[a)]
\item
If $Y$ is $T_1$ and $f:X\to Y$ is a connectedness preserving
function then $f$ is continuous.
\item
If $Y$ is $T_1$ and $f:X\to Y$ is a preserving function then
$f$ is continuous (i.e.\ $Pr(X,T_1)$ holds).
\item
If $Y$ has the cofinite topology and $f:X\to Y$ is a preserving
function then $f$ is continuous.
\item
If $A\subset X$ is not closed then there exists a connected set
$H\subset X$ such that $H\cap A \not = \emptyset \not = H - A$ and
$H-A$ is finite.
\end{enumerate}
\end{theorem}

\begin{proof} a) $\Rightarrow$ b) and b)$\Rightarrow$
c) are obvious.

c) $\Rightarrow$ d) Assume that $A\subset X$ is not closed. Let
$Y$ denote the space with the cofinite topology on the underlying
set of $X$. Choose a point $a_0 \in A$. ($A$ is not closed so it
is not empty, either.) Define now the function $f:X\to Y$ by
\begin{displaymath}
f(x)=\left\{
\begin{array}{ll}
a_0&\text{ if $x\in A$,}\\
x, &\text{otherwise.}
\end{array}
\right.
\end{displaymath}
Then $f$ is not continuous because the inverse image of the closed set
$\{a_0\}$ is the non-closed set $A$ hence ,by c), $f$ is not
preserving. As an arbitrary subset of $Y$ is compact, $f$
preserves compactness, so there is a connected set $H\subset X$
such that $f(H)$ is not connected. It follows that $H$ is infinite
and $f(H)$ is finite but not a singleton. As $f$ is the identity map on
$X-A$, the set $H-A$ is finite and so $H\cap A \not =
\emptyset$.Finally, $H\subset A$ is impossible because $f(H)$ is
not a singleton.

d) $\Rightarrow$ a) Assume $f:X\to Y$ is not continuous for a
$T_1$-space $Y$, hence there is a closed set $F\subset X$ such
that $A=f^{-1}(F)$ is not closed in $X$. By d), there is a
connected set $H$ such that $H\cap A \not = \emptyset$ and
$\emptyset \not = H - A$ is finite. But then $f(H)$ is not
connected because it is the the disjoint union of two non-empty
relatively closed sets, namely of $f(H)\cap F$ and of the finite set
$f(H)-F$. Consequently, $f$ does not preserve connectedness.
\end{proof}

\begin{corollary}
\label{5.2}
If $Pr(X,T_1)$ holds for a $T_1$-space $X$ then every closed subspace of
$X$ is the topological sum of its components.
\end{corollary}

\begin{proof} 
Let $K$ be a component of the closed subset $F\subset X$. 
It is enough to prove that $K$ is relatively open in $F$. 
Assume this is false; then $A=F-K$ is not closed in $X$, and thus, by
condition d) of Theorem \ref{5.1}, there is a connected set $H$ in $X$
such that $H\cap A \not = \emptyset$ and $\emptyset \not = H-A$ is
finite. 
Then $H - F$ is also finite, consequently $H\subset F$ because $H$ is
connected and $F$ is closed. 
Thus $H$ is a connected subset of $F$ which meets the component $K$ of
$F$, contradicting that $H\cap A \not = \emptyset$.
\end{proof}

\begin{corollary}
\label{5.3} 
If $Pr(X,T_1)$ holds for a $T_3$-space $X$ then every closed subspace of
$X$ is locally connected.
\end{corollary}

\begin{proof} 
By \ref{5.2} it is enough to prove that if every closed subset of a
regular space $X$ is the topological sum of its components then $X$ is
locally connected.

Let $U$ be a closed neighbourhood of a point $x\in X$.
By our assumption if $K$ denotes the component of $x$ in $U$ then $K$ is
open in $U$, hence $K\subset U$ is a connected neighbourhood of $x$ in 
$X$. 
As the closed neighbourhoods of a point form a neighbourhood-base of
the point in a regular space, $X$ is locally connected.
\end{proof}

\begin{theorem}
\label{5.4} 
If $Pr(X,T_1)$ holds for a $T_3$-space $X$ then $X$ is discrete.
\end{theorem}

By Corollary \ref{5.3} it is enough to prove the following result that,
we think, is interesting in itself:

\begin{theorem}
\label{5.5} 
If $X$ is $T_3$ and every regular closed subspace of $X$ is locally
connected then $X$ is discrete.
\end{theorem}

\begin{proof} 
We can assume without any loss of generality that $X$ is 
connected. Suppose, indirectly, that $x$ is a non-isolated point in $X$. 
By regularity, there is a sequence of non-empty open sets 
$\{G_n: n<\omega \}$ such that $x \not \in G_n$ and 
$\overline G_n \subset G_{n+1}$ for all $n<\omega$. 
Then the open set $G= \bigcup \{G_n : n < \omega \}$ can not be also
closed in the connected space $X$, so there is a point 
$p \in \overline G - G$.

Put $U_0=G_0$ and $U_n=G_n - \overline G_{n-1}$ for $0<n<\omega$.
If $H_0= \bigcup \{U_n : n \hbox { is even}\}$ and $H_1= \bigcup
\{U_n : n \hbox { is odd}\}$,
we claim that then $\overline G = \overline {H_0
\cup H_1} = \overline {H_0} \cup \overline {H_1} $.
Indeed, if $z\in \overline G$ and $W$ is any open
neighbourhood $W$ of $z$ then there is a
$k<\omega$ such that $W\cap G_k \not = \emptyset$.
Clearly, if $n$ is the
least such $k$ then we have $\emptyset \not = W\cap U_n$
as well, and so $x\in \overline {H_0 \cup H_1}$.

Consequently, $p \in \overline H_0$ or $p \in \overline H_1$;
assume e. g. that $p \in \overline H_0$. We shall show that then
$\overline H_0$ is not locally connected at $p$, although it
is a regular closed set, arriving at a contradiction.

Indeed, let $U$ be any neighbourhood of $p$ in $\overline{H_0}$
and fix an even number $n<{\omega}$ with $U_n\cap U=\emptyset$.
Then $U\cap U_{n+1}\subset \overline{H_0}\cap H_1=\emptyset$
implies $U\subset \overline{G_n}\cup (X\setminus G_{n+1})$,
where $\overline{G_n}$ and $X\setminus G_{n+1}$ are disjoint closed
sets both meeting $U$, hence $U$ is disconnected.
\end{proof}

With a little more effort it can also be shown that for any non-isolated
point $p$ in a $T_3$ space $X$ there is a regular closed set $H$ in $X$
with $p \in H$ such that $p$ is not a local connectivity point in $H$.

We do not know if every $T_2$ space $X$ with $Pr(X,T_1)$ has to be 
discrete. Also, the following $T_2$ version of Theorem \ref{5.5} seems to 
be open:
{\em If $X$ is $T_2$ and all closed subspaces of $X$ are locally connected
then $X$ has to be discrete.}
Note that if $X$ has the cofinite topology then it is hereditarily locally
connected and satisfies $Pr(X,T_1)$.

%\bibliographystyle{amsplain}
%\bibliography{10}
\providecommand{\bysame}{\leavevmode\hbox to3em{\hrulefill}\thinspace}
\providecommand{\MR}{\relax\ifhmode\unskip\space\fi MR }
\providecommand{\MRhref}[2]{%
  \href{http://www.ams.org/mathscinet-getitem?mr=#1}{#2}
}
\providecommand{\href}[2]{#2}

\end{document}